\documentclass{elsart}
\usepackage{amssymb,amscd}
\usepackage{latexsym}
\usepackage{euscript}
\usepackage{pb-diagram}

\newtheorem{e-proposition}[theorem]{Proposition}

\newtheorem{e-definition}[theorem]{Definition\rm}

\begin{document}
\begin{frontmatter}


\author{}
\ead{bzhao@math.wustl.edu}
\title{Dirichlet forms on hyperfinite $II_1$ factor}
\author{Bo Zhao}
\address{Department of Mathematics, Washington University in St. Louis, St. Louis,MO 63130}
\begin{abstract} Based on the structure of the hyperfinite $II_1$
factor, we study its Dirichlet forms which can be constructed from
Dirichlet forms on $M_{2^n}(\mathbb{C})$.

\end{abstract}

\end{frontmatter}
------------------------------------------------------------------------------

\section*{Introduction}
The aim of this paper is to study noncommutative Dirichlet forms on
hyperfinite $II_1$ factor. Dirichlet form in the commutative
setting, originated with the work of Beurling and Deny\cite{AB}, has
been developed especially by Fukushima\cite{MF} and
Silverstein\cite{MS}. The corresponding noncommutative theory, in
the $C^*$-algebra setting, is introduced by Albeverio and H\o
egh-Krohn\cite{SR}. It has been recognized that it shares a flavor
of geometry in the sense of Connes' noncommutative
geometry\cite{AC}. For a recent account of the theory, we refer the
reader to\cite{CF}\cite{CS}.

By the uniqueness of the hyperfinite $II_1$ factor $\mathcal{R}$, we
regard it as the completion of $\cup_n M_{2^n}(\mathbb{C})$ under
the normalized tracial state $\tau$. By the virtue of Gelfand and
Naimark theorem, each commutative $C^*$ algebra is isometrically
*-isomorphic to the algebra of continuous functions $C(X)$. Moreover,
each commutative von Neumann algebra is isometrically *-isometric to
the algebra of bounded functions $L^{\infty}(X)$. Combining the
above facts together, we make the following correspondence.

\begin{tabular}{|l|r|}
  \hline
   Noncommutative & Commutative \\
  $\cup_j M_{2^j}(\mathbb{C})$ & $\cup_j C(\mathcal{C}_j)$ \\
  CAR algebra & $C(\mathcal{C})$ \\
  $\mathcal{R}$ & $L^{\infty}(\mathcal{C})$ \\
  $L^2(\mathcal{R},\tau)$ & $L^{2}(\mathcal{C})$ \\
  \hline
\end{tabular}

where the notation $\mathcal{C}_j$ means $j$-th step Cantor set,
namely, $\mathcal{C}_j=\{\sum_{i=1}^j a_i 3^{-i}| a_i=0,2\}$ and
$\mathcal{C}$ means the cantor set.

Cantor set is a fractal. From that point of view, it can be obtained
from the self-similar sets. The Dirichlet form as well as harmonic
analysis on fractals has already been studied by Kigami\cite{JK},
Kusuoka\cite{SK}. Its basic idea is to construct a sequence of
Dirichlet forms on an increasing sequence of finite sets in such a
way that they satisfy a certain compatibility condition.
Consequently, we could obtain a Dirichlet form on the closure of
this increasing finite sets with some metric. Unfortunately, in
general, this closure is merely a proper subset of this fractal in
certain cases.

From the above correspondence, to study Dirichlet forms on
$L^2(\mathcal{R},\tau)$, we can first examine Dirichlet forms on
$M_{2^j}(\mathbb{C})$. However, positive operators on
$M_{2^j}(\mathbb{C})$ are still unclear. All that we know at present
is a nice machinery to write down systematically many interesting
positive maps. But, the completely positive map on
$M_{2^j}(\mathbb{C})$ is well understood. It has the form
$$\triangle_j(a)= \sum_{i=1}^{k} [[m_i,[m_i, a]]+ha+ah,$$ where
$m_i=m_i^*\in M_{2^j}$ and $h=h^*\in M_{2^j}$. The Dirichlet form
associated with a uniformly continuous completely positive semigroup
on a type $I$ von Neumann algebra was completely understood
\cite{AG}.

The paper is organized as follows. In section 1, we give a brief
exposition of the noncommutative Dirichlet forms theory. Section 2
establishes the relation between the Dirichlet forms on
$\mathcal{R}$ and Dirichlet forms on $M_{2^j}(\mathbb{C})$. It will
be shown that if $\mathcal{E}$ is a bounded, then it is the limiting
of the Dirichlet form on $M_{2^j}(\mathbb{C})$. And if $\cup_j
M_{2^j}(\mathbb{C})$ is the form core, then $\mathcal{E}$ can be
recovered from Dirichlet forms on $M_{2^j}(\mathbb{C})$. In section
3, we will give a concrete example.

I wish to thank  Nik Weaver. Without his many helpful suggestions
this work could not have been done.

\section {Dirichlet forms} \label{diri}
 Let $\mathcal{A}$ be a $C^*$-algebra and
$\mathcal{A}^{**}$ its enveloping von Neumann algebra, with unit
$1_{\mathcal{A}^{**}}$. Let's consider a densely defined, faithful,
semifinite, lower semicontinuous trace $\tau$ on $\mathcal{A}$. We
denote by $L^2(\mathcal{A},\tau)$
($\langle,\rangle_{L^2(\mathcal{A},\tau)}$) the Hilbert space of the
GNS representation $\pi_{\tau}$ associated to $\tau$, and by
$L^\infty(\mathcal{A},\tau)$ or $\mathcal{M}$ the von Neumann
algebra $\pi_{\tau}(\mathcal{A})''$ in
$\mathcal{B}(L^2(\mathcal{A},\tau))$ generated by $\mathcal{A}$ in
the GNS representation. When unnecessary, we shall not distinguish
between $\tau$ and its canonical normal extention on $\mathcal{M}$,
between elements on $\mathcal{A}$ and their representation in
$\mathcal{M}$ as a bounded operator on $L^2(\mathcal{A},\tau)$, nor
between elements $a$ of $\mathcal{A}$ or $\mathcal{M}$ which are
square integrable $(\tau(a*a)<+\infty)$ and their canonical image in
$L^2(\mathcal{A},\tau)$. Then $||a||$ stands for the uniform norm of
$a$ in $\mathcal{A}$ or in $\mathcal{M}$, $||a||_2$ or
$||a||_{L^2(\mathcal{A},\tau)}$ for $L^2$-norm of $a$ in
$L^2(\mathcal{A},\tau)$, $1_{\mathcal{M}}$ for the unit of
$\mathcal{M}$. As usual, $L^{\infty}_+(\mathcal{A},\tau)$ and
$L^{2}_+(\mathcal{A},\tau)$ will denote the positive part of
$L^{\infty}(\mathcal{A},\tau)$ and $L^{2}(\mathcal{A},\tau)$,
respectively.

Recall that
$(\mathcal{M},L^{2}(\mathcal{A},\tau),L^{2}_+(\mathcal{A},\tau))$
is a \emph{standard form} of the von Neumann algebra
$\mathcal{M}$. In particular $L^{2}_+(\mathcal{A},\tau)$ is a
closed convex cone in $L^{2}(\mathcal{A},\tau)$, inducing an
anti-linear isometry (the modular conjugation) $J$ on
$L^{2}(\mathcal{A},\tau)$ which is the extension of the involution
$a\mapsto a^*$ of $\mathcal{M}$. The subspace of $J$-invariant
elements (called \emph{real}) will be denoted by
$L^{2}_h(\mathcal{A},\tau)$.

When $a$ is real, the symbol $a\wedge 1$ will denote the Hilbert
projection onto the closed and convex subset $C$ of
$L^{2}_h(\mathcal{A},\tau)$ obtained as the $L^2$-closure of
$\{a\in L^{2}_+(\mathcal{A},\tau): a\leq 1_{\mathcal{M}}\}$.

\newtheorem{guess}{Definition}[section]
\begin{guess}\label{sg}
Given a strongly continuous semigroup $\Phi_t (t\in R^+)$ of
operators defined on $L^\infty(\mathcal{A},\tau)$,
\begin{enumerate}
    \item it is
\emph{symmetric}, if $\tau( \Phi_t(x)y)=\tau(x \Phi_t(y))$.
    \item it is \emph{Markov}, if $0\leq x\leq 1_{\mathcal{M}}$ implies that $0\leq
\Phi_t(x)\leq 1_{\mathcal{M}}$.
    \item it is \emph{conservative}, if
$\Phi_t(1_{\mathcal{M}})=1_{\mathcal{M}}$.
    \item it is \emph{completely positive}, if for any $n$ we have $\sum_{i,j=1}^n
b_i^* \Phi_t(a_i^* a_j)b_j\geq 0$ where $a_i,b_i\in \mathcal{M}
,i=1,\ldots,n$.
\end{enumerate}
\end{guess}

\newtheorem{form}[guess]{Definition}
\begin{form}\label{defi}
A closed, densely defined, nonnegative quadratic form
$(\mathcal{E},D(\mathcal{E}))$ on $L^2(\mathcal{A},\tau)$ is said
to be
\begin{enumerate}
    \item \emph{real} if for $a\in D(\mathcal{E})$,$J(a)\in
    D(\mathcal{E})$ and $\mathcal{E}(J(a))=\mathcal{E}(a)$;
    \item a \emph{Dirichlet form} if it is real and
     \[\mathcal{E}(a\wedge 1)\leq \mathcal{E}(a),\textit{ for }
     a\in D(\mathcal{E})\cap L^2_h(\mathcal{A},\tau);\]
    \item a \emph{completely Dirichlet form} if the canonical extension
    $(\mathcal{E}^n, D(\mathcal{E}^n))$ to $L^2 (M_n(\mathcal{A}),\tau_n)$
     \[\mathcal{E}^n[[a_{i,j}]_{i,j=1}^n]:=\sum_{i,j=1}^n\mathcal{E}[a_{i,j}]\emph{, where }[a_{i,j}]_{i,j=1}^n\in D(\mathcal{E}^n)):=M_n(D(\mathcal{E}))\]
is a Dirichlet form for all $n\geq 1$.
\end{enumerate}
\end{form}

By the general theory, a symmetric, Markov semigroup gives rise to a
Dirichlet form $\mathcal{E}$ on $L^2(\mathcal{A},\tau)$ by
\[\mathcal{E}(a)=\langle \triangle a, a\rangle_2,\] where $\triangle$
is the generator of this semigroup. And completely positive,
symmetric, Markov semigroup gives rise to a completely Dirichlet
form. On the other hand, starting from a Dirichlet form on
$L^2(\mathcal{A},\tau)$, one can always reconstruct a symmetric,
Markov semigroup. And furthermore, if this Dirichlet form is
complete, then this semigroup is completely positive.

To summarize, conditions (1)(2) from Definition \ref{sg} is
equivalent to (1)(2) from Definition \ref{defi}, and (1)(2)(4) from
Definition \ref{sg} is equivalent to (1)(2)(3) from Definition
\ref{defi}.

Given $\mathcal{E}$, define the inner product $\langle, \rangle_1$
on $D(\mathcal{E})$ by
\[\langle a, b\rangle_1=\langle \triangle a, b\rangle_2+ \langle a, b\rangle_2.\]
The form $\mathcal{E}$ is closed implies $\langle, \rangle_1$ is a
Hilbert space.

\section {Dirichlet forms on hyperfinite $II_1$ factor} \label{matrix}
If $\mathcal{R}$ is the hyperfinite $II_1$ factor, from the property
that the hyperfinite $II_1$ factor is uniquely determined up to
 isomorphism, we regard $\mathcal{R}$ as the completion of $\cup_n M_{2^n}(\mathbb{C})$ under
the normalized tracial state $\tau$. The mapping
 \[a\rightarrow \left[
\begin{array}{rr}
 a & 0   \\
                  0 & a   \\
\end{array}
\right] \]
 is the embedding of $M_{2^n}(\mathbb{C})$ as a subalgebra of
                    $M_{2^{n+1}}(\mathbb{C})$. Throughout the paper, the notation
 $\tau_n$
                    means the normalized tracial state on
                    $M_{2^n}(\mathbb{C})$.

\newtheorem{same}[guess]{Lemma}
\begin{same}
$\cup_n M_{2^n}(\mathbb{C})$ is dense in $L^{2}(\mathcal{R},\tau)$.
\end{same}
\begin{pf}
Given $a\in L^2$, we find a sequences $a_n\in \mathcal{R}$, such
that $a_n\rightarrow a$ in $L^2$, for each $a_n$, we find $b_n\in
M_{2^n}(\mathbb{C})$, such that $\langle
(b_n-a_n)x,x\rangle_{L^2}<\frac{1}{n}, \forall x\in L^2$. Then the
sequences $b_n-a_n\in \mathcal{R}$ weak operator converges to $0$,
and so it's strong operator converges to $0$. In particular,
$||b_n-a_n||_{L^2}\rightarrow 0$. Then
$||b_n-a||_{L^2}=||b_n-a_n+a_n-a||_{L^2}\leq
||b_n-a_n||_{L^2}+||a_n-a||_{L^2}\rightarrow 0$ when $n\rightarrow
\infty$.
\end{pf}

For each $n$, we have a conditional expectation map
\[E_n: L^2\rightarrow M_{2^n}(\mathbb{C})\] and an extension map
\[\Pi_n: M_{2^n}(\mathbb{C})\rightarrow L^2.\]
On $\cup_{j\geq 0} M_{2^j}(\mathbb{C})$, the map $E_n$  is defined
by \[E_n(a_1\otimes a_2\otimes \ldots \otimes a_n\otimes
\ldots\otimes a_{n+k})=(\prod_{j=1}^k \tau_1(a_{n+j})) a_1\otimes
a_2\ldots \otimes a_n,\] since $\cup_{j\geq 0} M_{2^j}(\mathbb{C})$
is dense in $L^2$, we then extend $E_n$ to $L^2$. The map $\Pi_n$ is
defined from the above embedding map.

 Let $P_n=\Pi_n\circ E_n$, $Q_n=I-P_n$, then
it's not hard to see $P_n$ (resp., $Q_n$) is a family of projections
which increase (resp., decrease) to $I$ (resp., $0$) when
$n\rightarrow \infty$.

Given a Dirichlet form $\mathcal{E}$ on $L^2$, unless otherwise
stated, we assume that $D(\mathcal{E})\supset \cup_{j\geq 0}
M_{2^j}(\mathbb{C})$. The point of using this assumption is that the
dense subalgebra $\cup_{j\geq 0} M_{2^j}(\mathbb{C})$ is realized as
the smooth functions from the commutative point of view.

\newtheorem{cpd}[guess]{Proposition}
\begin{cpd}
Let $\mathcal{E}_n(a):=\mathcal{E}(P_n a)$, then $\mathcal{E}_n$ is
a bounded Dirichlet form on $L^2$.
\end{cpd}
\begin{pf}
The proof is straightforward.
\end{pf}

\newtheorem{sam}[guess]{Proposition}
\begin{sam}
Given a Dirichlet form $\mathcal{E}$ on $L^2$, the following
conditions are equivalent:
\begin{enumerate}
  \item $\lim _{n\rightarrow
\infty}\mathcal{E}_n(a)=\mathcal{E}(a),\forall a\in D(\mathcal{E})$.
  \item $\lim _{n\rightarrow
\infty}\mathcal{E}(Q_n a)=0,\forall a\in D(\mathcal{E})$
\end{enumerate}\label{equi}
 \end{sam}
\begin{pf}
For the given $\mathcal{E}$, let $\triangle$ be its Markovian
semigroup generator.\\
 (2)$\Rightarrow$ (1):
From the triangle inequality, we have
\begin{eqnarray*}
\mathcal{E}^{\frac{1}{2}}(Q_na)&=& ||\triangle^{\frac{1}{2}}Q_n
a||_2=||\triangle^{\frac{1}{2}} a-\triangle^{\frac{1}{2}}P_n a||_2
\geq \left| ||\triangle^{\frac{1}{2}}P_n
a||_2-||\triangle^{\frac{1}{2}} a||_2 \right |
\\ &\geq& | \mathcal{E}^{\frac{1}{2}}_n(a)- \mathcal{E}^{\frac{1}{2}}(a)
|.
\end{eqnarray*}
 By letting $n\rightarrow
\infty$ to above inequality gives (1).

(1)$\Rightarrow$ (2): From the spectral representation of
$\triangle$, we have
\[\mathcal{E}(a)=\langle\triangle a, a \rangle=\int_0^{\infty}\lambda d\langle F_{\lambda} a, a \rangle,\]
where $F_{\lambda}$ is the spectral projection of $\triangle$. Write
$f_n(\lambda)=\langle F_{\lambda}P_n a,P_n a \rangle$,
$g_n(\lambda)=\langle F_{\lambda}Q_n a,Q_n a \rangle$,
$f(\lambda)=\langle F_{\lambda}a,a \rangle$. This gives
\[\mathcal{E}_n(a)=\langle\triangle P_n a,P_n a \rangle=\int_0^{\infty}\lambda d\langle F_{\lambda}P_n a,P_n a \rangle=\int_0^{\infty}\lambda df_n(\lambda).\]
\[\mathcal{E}(Q_na)=\langle\triangle Q_n a,Q_n a \rangle=\int_0^{\infty}\lambda d\langle F_{\lambda}Q_n a,Q_n a \rangle=\int_0^{\infty}\lambda dg_n(\lambda).\]
 The functions $f_n,
g_n, f$  are increasing imply they are differentiable a.e. Denote
\[f'(\lambda)=\lim_{h\rightarrow
0}\frac{f(\lambda+h)-f(\lambda)}{h}=\langle\lim_{h\rightarrow 0}
\frac{F_{\lambda+h}-F_{\lambda}}{h}a,a \rangle=\langle
F'_{\lambda}a,a \rangle\] where $F'_{\lambda}$ is an unbounded
positive operator on $L^2$. Next we observe that
\[\sqrt{g'_n(\lambda)}=||\sqrt{ F'_{\lambda}}Q_n a||_2=||\sqrt{ F'_{\lambda}}a-\sqrt{ F'_{\lambda}}P_na||_2\leq ||\sqrt{ F'_{\lambda}}a||_2+||\sqrt{ F'_{\lambda}}P_na||_2,\]
which gives
\[\lambda g'_n(\lambda)=\lambda \langle F'_{\lambda}Q_n a, Q_n a\rangle \leq 2\lambda \langle  F'_{\lambda}a, a\rangle+2\lambda \langle F'_{\lambda}P_na, P_na\rangle=2\lambda f'
(\lambda)+2\lambda f'_n(\lambda).\] Notice $f_n\rightarrow f$,
$g_n\rightarrow 0$ as $n\rightarrow \infty$, this gives $\lambda
f'_n(\lambda)\rightarrow \lambda f'(\lambda)$, $\lambda
g'_n(\lambda)\rightarrow 0$ a.e. when $n\rightarrow \infty$. To
complete the proof, we apply the generalized dominated convergence
theorem. Thus
\[\lim_{n\rightarrow \infty}\mathcal{E}(Q_na)=\lim_{n\rightarrow \infty}\int_0^{\infty}\lambda g'_n(\lambda)d\lambda =\int_0^{\infty}\lim_{n\rightarrow \infty}\lambda g'_n(\lambda)d\lambda=0\]
\end{pf}

The remainder of this section will be devoted to discuss the
relation between $\mathcal{E}$ and $\mathcal{E}_n$.

\newtheorem{bd}[guess]{Theorem}
\begin{bd}
If $\mathcal{E}$ is a bounded Dirichlet form on $L^2$, then
\[\lim_{n\rightarrow \infty}\mathcal{E}_n(a)=\mathcal{E}(a),\forall a\in
L^2.\]
\end{bd}
\begin{pf}
The Markov semigroup generator $\triangle$ of $\mathcal{E}$ is
bounded. We have
\[\mathcal{E}(Q_na)=||\sqrt{\triangle}Q_na||_2^2\leq ||\triangle||||Q_na||_2^2\rightarrow 0, \textrm{when } n\rightarrow \infty,\]
which completes the proof from proposition \ref{equi}.
\end{pf}

\newtheorem{fc}[guess]{Theorem}
\begin{fc}
Given a Dirichlet form on $L^2$, if $\lim _{n\rightarrow
\infty}\mathcal{E}_n(a)=\mathcal{E}(a),\forall a\in D(\mathcal{E})$,
then $\cup_{j\geq 0} M_{2^j}(\mathbb{C})$ is the form core.
\end{fc}
\begin{pf}
For $a\in D(\mathcal{E})$, $\lim _{n\rightarrow
\infty}\mathcal{E}_n(a)=\mathcal{E}(a)$ gives $\lim _{n\rightarrow
\infty}\mathcal{E}(Q_na)=0$ from proposition \ref{equi}. It implies
$P_na\in M_{2^n}(\mathbb{C})\rightarrow a$ in $\langle.,.\rangle_1$
norm, which completes the proof.
\end{pf}

For the other direction of the above theorem, we have the following
result.

\newtheorem{fd}[guess]{Theorem}
\begin{fd}
Given a Dirichlet form $\mathcal{E}$ on $L^2$, if $\cup_{j\geq 0}
M_{2^j}(\mathbb{C})$ is the form core, then $\mathcal{E}$ can be
recovered from $\mathcal{E}_n$.
\end{fd}
\begin{pf}
Given $\mathcal{E}$, we first construct $\mathcal{E}_n$, and then
define
\[\mathcal{F}(a)=\lim_{n\rightarrow\infty}\mathcal{E}_n(a), \forall a\in \cup_{j\geq 0}
M_{2^j}(\mathbb{C}).\] It is well defined, because eventually
$\mathcal{E}_n(a)$ will stay the same when $n$ is large enough.
Since $\mathcal{E}$ and $\mathcal{F}$ agree on $\cup_{j\geq 0}
M_{2^j}(\mathbb{C})$ and  $\cup_{j\geq 0} M_{2^j}(\mathbb{C})$ is
the form core, after the completion of $\cup_{j\geq 0}
M_{2^j}(\mathbb{C})$ under $\langle, \rangle_1$ norm, we get
$\mathcal{F}$ and consequently $\mathcal{E}$, which gives the proof.
\end{pf}

The principal significance of the above theorem is that it allows
one to construct a certain type Dirichlet forms on $L^2$. The
procedure is: first we construct $\tilde{\mathcal{E}}_n$ on
$M_{2^n}(\mathbb{C})$ in such a way that $\tilde{\mathcal{E}}_n$ is
compatible with $\tilde{\mathcal{E}}_{n+1}$ on
$M_{2^n}(\mathbb{C})$, i.e.,
$\tilde{\mathcal{E}}_n(a)=\tilde{\mathcal{E}}_{n+1}(a), \forall a\in
M_{2^n}(\mathbb{C})$. Next, we define
$\mathcal{E}_n(a)=\tilde{\mathcal{E}}_n(E_na)$, so that
$\mathcal{E}_n$ is defined on $L^2$. Then, we
 define
\[\mathcal{E}(a)=\lim_{n\rightarrow\infty}\mathcal{E}_n(a), \forall a\in \cup_{j\geq 0}
M_{2^j}(\mathbb{C}).\] Finally, we take the completion of
$\cup_{j\geq 0} M_{2^j}(\mathbb{C})$ under $\langle, \rangle_1$
norm. The constructed Dirichlet form $\mathcal{E}$ has the property
that $\cup_{j\geq 0} M_{2^j}(\mathbb{C})$ is its form core.

\section {An example}
 As it is known that $L^{\infty}[0,1]$ is the maximal
abelian subalgebra of $\mathcal{R}$, thus\cite{MT}, we have a
conditional expectation map
\[B: \mathcal{R}\rightarrow L^{\infty}[0,1].\]
Intuitively, this conditional expectation map is coming from the
following diagram.
\[
\begin{CD}
&M_{2^1} \subset  &M_{2^2} \subset  & M_{2^3} \subset &\ldots
\subset
& \mathcal{R}\\
&\downarrow{B_1} &\downarrow{B_2}
&\downarrow{B_3}&&\downarrow{B}\\
&D_{2^1} \subset  &D_{2^2} \subset  & D_{2^3} \subset &\ldots
\subset
& L^{\infty}[0,1]\\
\end{CD}
\]
where $D_{2^n}$ is the $2^n\times 2^n$ diagonal matrix and $B_n$ is
the conditional expectation map from $M_{2^n}$ to $D_{2^n}$.

For this conditional expectation $B$, we can extend it to
$L^2(\mathcal{R},\tau)$ since $\mathcal{R}$ is dense in
$L^2(\mathcal{R},\tau)$. Keep the same notation, let this map be $B:
L^2(\mathcal{R},\tau)\rightarrow L^2[0,1] $.
\newtheorem{di}[guess]{Proposition}
\begin{di}
For the conditional expectation operator $B$, we define
\[\mathcal{E}(a)=\langle (I-B)a, a\rangle_2,\] then $\mathcal{E}$ is a bounded, completely Dirichlet form on $L^2$.
\end{di}
\begin{pf}
Let $\triangle=I-B$, then it's a projection operator, hence bounded.
Now it's obvious that the semigroup $\phi_t=e^{-t\triangle}$ is
Markov and symmetric. To show it's a completely Dirichlet form,
notice from Definition \ref{defi}, the generator of $\mathcal{E}^n$
is $\triangle\otimes I_n(\mathbb{C})$ which is again a projection.
This completes the proof.
\end{pf}

From section (2), for this Dirichlet form $\mathcal{E}$, we get the
restricted $\mathcal{E}_n$. And it's not hard to see that
\[\mathcal{E}_n(a)=\sum _{i=1}^{2^n}\tau_n([p_i, E_n a][p_i, E_na]^*).\]
where $p_i$ is the diagonal matrix with $1$ in entry $(i,i)$ and $0$
elsewhere. The generator $\triangle_n$ of $\mathcal{E}_n$ is
$\triangle_n(a)=\sum_{j=1}^{2^n}[p_j,[p_j,E_na]]$. As shown from
section 2, $\triangle_n$ weak operator converges to $\triangle$.

From \cite{S}\cite{CS}, for this Dirichlet form $\mathcal{E}$, it is
naturally equipped with a structure of a Hilbert bimodule over
$\mathcal{R}$, and a derivation operator $\partial$. In order to
understand the bimodule structure, we begin with considering the
Dirichlet form $\mathcal{E}_n$. For this Dirichlet form
$\mathcal{E}_n$, the associated derivation map
\[\partial_n:\mathcal{R}\rightarrow  \mathcal{R}\otimes \mathbb{C}^{2^n}.\]
is defined by
\[\partial_n(a)=\oplus_{j=1}^{2^n}[p_j,E_na]\textrm{, for } a\in
\mathcal{R}.\] $\partial_n$ is used to define the noncommutative
differential calculus in \cite{22}. The Hilbert bimodule is
$\mathcal{R}\otimes \mathbb{C}^{2^n}=l^2(2^n,\mathcal{R})$ and its
structure is
\[(f\cdot a)(j)=f(j)a, (a\cdot f)(j)=af(j),\textrm{ for } f\in l^2(2^n,\mathcal{R}),\]
\[\langle f, g \rangle_{\mathcal{R}}=\sum_{j=1}^{2^n} f(j)^*g(j),\textrm{ for } f,g\in l^2(2^n,\mathcal{R}). \]

Return to $\mathcal{E}$, it follows in the same manner, the Hilbert
bimodule associated with $\mathcal{E}$ is $\mathcal{R}\otimes
L^2[0,1]=L^2([0,1],\mathcal{R})$ and its bimodule structure is
\[(f\cdot a)(t)=f(t)a, (a\cdot f)(t)=af(t),\textrm{ for } f\in L^2([0,1],\mathcal{R}),\]
\[\langle f, g \rangle_{\mathcal{R}}=\int_0^1 f(t)^*g(t)dt,\textrm{ for } f,g\in L^2([0,1],\mathcal{R}). \]
The derivation map is
\[\partial:\mathcal{R}\rightarrow
\mathcal{R}\otimes L^2[0,1]=L^2([0,1],\mathcal{R}). \]

\newtheorem{de3}[guess]{Definition{\cite{S}}}
\begin{de3}
A completely positive Markov semigroup $(\Phi_t)_{t\geq 0}$ and
its infinitesimal generator $\triangle$ are \emph{strongly local}
if the associated Hilbert bimodule is trivialisable, that is
isometrically imbeded in an amplication of the trivial bimodule.
\end{de3}

\newtheorem{local}[guess]{Proposition}
\begin{local}
The Dirichlet form $\mathcal{E}$ defined above is strongly local.
\end{local}
\begin{pf}
This is because the bimodule associated to $\mathcal{E}$ is
$\mathcal{R}\otimes L^2[0,1]$.
\end{pf}

\bibliographystyle{amsplain}

\end{document}